\documentclass[11pt,a4paper]{article}
\usepackage[USenglish]{babel}
\usepackage[T1]{fontenc}
\usepackage[utf8]{inputenc}
\usepackage{csquotes}
\usepackage{enumitem}
\usepackage{appendix}
\usepackage[colorlinks=true,linkcolor=blue,bookmarksnumbered=true,bookmarksdepth=3,pdftitle=Representation\ Formula\ for\ Viscosity\ Solutions\ to\ Parabolic\ PDEs\ with\ Sublinear\ Operators,pdfauthor=Marco\ Pozza]{hyperref}

\usepackage{sectsty}
\subsubsectionfont{\normalfont\itshape}

\usepackage{amssymb}
\usepackage{amsmath}
\usepackage{mathtools}
\usepackage{amsthm}
\usepackage{mathrsfs}
\usepackage{dsfont}

\usepackage[noabbrev]{cleveref}

\theoremstyle{plain}
\newtheorem{theo}{Theorem}[section]
\crefname{theo}{theorem}{theorems}
\crefname{dpp}{dynamic programming principle}{}
\newtheorem{prop}[theo]{Proposition}

\newtheorem{lem}[theo]{Lemma}
\theoremstyle{definition}
\newtheorem{defin}[theo]{Definition}

\newtheorem{prob}[theo]{Problem}
\newtheorem{assum}[theo]{Assumptions}
\crefname{assum}{assumptions}{assumptions}
\theoremstyle{remark}
\newtheorem{rem}[theo]{Remark}

\DeclareMathOperator{\tr}{tr}
\DeclareMathOperator{\Lip}{Lip}
\DeclareMathOperator*{\esss}{ess\,sup}
\DeclareMathOperator*{\essi}{ess\,inf}
\DeclareMathOperator*{\argmin}{arg\,min}

\title{Representation Formula for Viscosity Solutions to Parabolic PDEs with Sublinear Operators}
\author{Marco Pozza\thanks{Dipartimento di Matematica ``G. Castelnuovo'', Sapienza Università di Roma. Piazzale Aldo Moro 5, 00185 Roma, Italy. {\tt pozza@mat.uniroma1.it}}}
\date{}

\begin{document}
\maketitle

\begin{abstract}
We provide a representation formula for viscosity solutions to a class of nonlinear second order parabolic PDE problem involving sublinear operators. This is done through a dynamic programming principle derived from \cite{art:pathpeng}. The formula can be seen as a nonlinear extension of the Feynman--Kac formula and is based on the backward stochastic differential equations theory.
\end{abstract}


\medskip


\section{Introduction}

It is well known that viscosity solutions were conceived by Crandall and Lions (1982) in the framework of optimal control theory. The goal was to show well posedness of Hamilton--Jacobi--Bellman equations in the whole space, and to prove, via dynamic programming principle, the value function of a suitable optimal control problem being the unique solution.

When trying to extend viscosity methods to the analysis of second order \emph{parabolic partial differential equations}, PDEs for short, and getting representation formulae, it appeared clear that some stochastic dynamics must be brought into play. Not surprisingly, this has been first done for stochastic control models. The hard work of a generation of mathematicians, \cite{hjlions2,hjlions1,bensoussa82,krylovbook,flemingrishel75,nisio15} among the others, allowed making effective dynamic programming approach to stochastic control problems.

Prompted by this body of investigations, a stream of research arose in the probabilistic community ultimately leading to the theory of \emph{backward stochastic differential equations}, BSDEs for short, which was introduced by Pardoux and Peng in \cite{Pardoux_Peng_1990} (1990). Since then, it has attracted a great interest due to its connections with mathematical finance and PDEs, as well as with stochastic control. This theory has been in particular used to extend the classical Feynman--Kac formula, which establishes a link between linear parabolic PDEs and \emph{stochastic differential equations}, SDEs for short, to semilinear and quasilinear equations, see for example \cite{Delarue_2002,Delarue_Guatteri_2006,Ma_Yong_2007}. See also \cite{pardouxbook} for a rather complete overview of the semilinear case.

For sake of clarity, let us consider the following semilinear parabolic PDE problem coupled with final conditions,
\begin{equation}\label{exparsemlin}
\begin{cases}
\begin{aligned}
\frac12\left\langle\sigma\sigma^\dag,D^2_xu\right\rangle(t,x)+(\nabla_xub)(t,x)&\\
+\partial_tu(t,x)+f(t,x,u,\nabla_xu\sigma)&=0
\end{aligned}
&t\in[0,T],x\in\mathds R^N,\\
u(T,x)=g(x),&x\in\mathds R^N,
\end{cases}
\end{equation}
then its viscosity solution can be written as $u(t,x)=\mathds E\left(Y^{t,x}_t\right)$, where $Y$ is given by the following system, called \emph{forward backward stochastic differential equation} or FBSDE in short, which is made in turn of two equations, the first one is a SDE, and the second one a BSDE depending on the first one
\[
\left\{\begin{aligned}
\!X^{t,x}_s\!=&x+\int_t^s\sigma\left(r,X^{t,x}_r\right)dW_r+\int_t^sb\left(r,X^{t,x}_r\right)dr,\\
\!Y^{t,x}_s\!=&g\left(X^{t,x}_T\right)\!+\!\!\int_s^T\!\!f\!\left(r,X_r^{t,x},Y_r^{t,x},Z_r^{t,x}\right)\!dr\!-\!\!\int_s^T\!\!Z_r^{t,x}dW_r,
\end{aligned}\right.s\!\in\![t,T],x\!\in\!\mathds R^N.
\]
As can be intuitively seen, the SDE takes care of the linear operator defined by $\sigma$ and $b$, also called the \emph{infinitesimal generators} of the SDE, while the BSDE depends on $f$ and $g$. In other words, this extension of the Feynman--Kac formula basically does not modify the treatment of the second order linear operator with respect to the completely linear case.

Subsequently, Peng introduced in \cite{math/0601035} (2006) the notion of \emph{$G$--expectation}, a nonlinear expectation generated by a fully nonlinear second order operator $G$ via its viscosity solutions. This work has originated an active research field, with relevant applications to Mathematical Finance.

Peng has improved this theory in several papers and has given a comprehensive account of it in the book \cite{book:peng}, where he highlights the role of the so--called \emph{sublinear expectations}, namely $G$--expectations generated by sublinear operators. Finally in \cite{art:pathpeng}, Peng provides representation formulae for viscosity solutions using these expectations. More precisely, given a sublinear operators $G$ and the $G$--heat equation
\begin{equation}\label{gheat}
\begin{cases}
\partial_tu(t,x)+G\left(D_x^2u(t,x)\right)=0,&t\in[0,T],x\in\mathds R^N,\\
u(T,x)=g(x),&x\in\mathds R^N.
\end{cases}
\end{equation}
he represents the viscosity solution as
\[
u(t,x)=\sup_{\sigma\in\mathcal A}\mathds E\left(g\left(x+\int_t^T\sigma_sdW_s\right)\right),
\]
where $\mathcal A$ is a family of stochastic process associated to $G$ and $W$ is a Brownian motion.
The key to prove it is a \emph{dynamic programming principle} that we will illustrate in the paper, see \cref{dynprincsec}. We point out that here the novelty with respect to the Feynman--Kac formula is essentially given by the sublinearity of the operator $G$.

The purpose of this article is to apply a generalized version of the dynamic programming principle of \cite{art:pathpeng} in order to give representation formulae of solutions to PDE problems of the type
\begin{equation}\label{parsubopgen}
\begin{cases}
\begin{aligned}
\partial_tu(t,x)+F\left(t,x,\nabla_xu,D^2_xu\right)&\\
+f(t,x,u,\nabla_xu)&=0,
\end{aligned}
&t\in[0,T],x\in\mathds R^N,\\
u(T,x)=g(x),&x\in\mathds R^N,
\end{cases}
\end{equation}
where $F$ is a sublinear operator, with respect the third and the fourth argument. This problem is clearly a blend between \eqref{exparsemlin} and \eqref{gheat}, where the additional difficulty with respect to \eqref{exparsemlin} is given by the sublinearity of the operator, while the generalization with respect to \eqref{gheat} is the dependence of $F$ on $(t,x)$ of and the presence of the term $f$.

This is hopefully just a first step to further extend the Feynman--Kac formula to problems with sublinear operators using a the BSDE theory in order to deal with general quasilinear problems. We also point out that there is close connection between this method and second order BSDEs, 2BSDEs for short. 2BSDEs were introduced by Cheridito, Soner, Touzi and Victoir in \cite{int2bsde} (2007). Then, in 2011, Soner, Touzi and Zhang \cite{sonertouzizhang11} provided a complete theory of existence and uniqueness for 2BSDEs under Lipschitz conditions. In those papers is also analyzed the connection between 2BSDEs and fully nonlinear PDEs. Among the subsequent developments of this theory we cite \cite{possamaitanzhou18,sonertouzizhang11,ekrentouzizhang16,ekrentouzizhang16_2,matoussipossamaizhou13,kazi-tanipossamaizhou15} and in particular \cite{possamai13}, which performs its analysis replacing the Lipschitz condition on $y$ with monotonicity as we do here.

This paper is organized as follows: in \cref{prelim} we make a preliminary study of the problem, analyzing the structure of the sublinear operator $F$ and developing a dynamic programming principle which is the core of our theory. Then, in \cref{parprob}, we perform the essential part of our analysis, and obtain in this way our main results. In \cref{2bsdesec} we summarily analyze the connection between 2BSDEs and our representation formula. The appendix at the end briefly gives some probability results we need, with a focus on the BSDE theory. We proceed setting the notation used in the paper.

\subsubsection*{Notation}
\pdfbookmark[3]{Notation}{Notation}

We will work on the filtered probability space $\left(\Omega,\mathcal F,\{\mathcal F_t\}_{t\in[0,\infty)},\mathds P\right)$,
\begin{itemize}
\item $\mathcal F$ is a complete $\sigma$--algebra on $\Omega$;
\item the stochastic process $\{W_t\}_{t\in[0,\infty)}$ will denote the $N$ dimensional Brownian motion under $\mathds P$;
\item $\{\mathcal F_t\}_{t\in[0,\infty)}$ is the filtration defined by $\{W_t\}_{t\in[0,\infty)}$ which respects the usual condition of completeness and right continuity;
\item $\{W_s^t\}_{s\in[t,\infty)}:=\{W_s-W_t\}_{s\in[t,\infty)}$ is a Brownian motion independent from $\{W_s\}_{s\in[0,t]}$ by the \emph{strong Markov property};
\item $\{\mathcal F^t_s\}_{s\in[t,\infty)}$ is the filtration generated by $\{W_s^t\}_{s\in[t,\infty)}$ which we assume respect the usual condition and is independent from $\mathcal F_t$;
\item we will say that a stochastic process $\{H_t\}_{t\in[0,\infty)}$ is adapted if $H_t$ is $\mathcal F_t$--measurable for any $t\in[0,\infty)$;
\item we will say that a stochastic process $\{H_t\}_{t\in[0,\infty)}$ is progressively measurable, or simply progressive, if, for any $T\in[0,\infty)$, the application that to any $(t,\omega)\in[0,T]\times\Omega$ associate $H_t(\omega)$ is measurable for the $\sigma$--algebra $\mathcal B([0,T])\times\mathcal F_T$;
\item a function on $\mathds R$ is called \emph{cadlag} if is right continuous and has left limit everywhere;
\item a cadlag (in time) process is progressive if and only if is adapted;
\item $B_\delta(x)$ will denote an open ball centered in $x$ with radius $\delta$;
\item for any Lipschitz continuous function $f$ we will denote its Lipschitz constant as $\Lip(f)$;
\item if $A\in\mathds R^{N\times M}$ then $A^\dag$ will denote its transpose and $\sigma_A$ its spectrum;
\item (Frobenius product) if $A,B\in\mathds R^{N\times M}$ then $\langle A,B\rangle:=\tr\left(AB^\dag\right)$ and $|A|$ is the norm $\sqrt{\langle A,A\rangle}=\sqrt{\sum\limits_{i=1}^N\sum\limits_{j=1}^MA_{i,j}^2}$;
\item $\mathds S^N$ is the space of all $N\times N$ real valued symmetric matrices and $\mathds S^N_+$ is the subset of $\mathds S^N$ made up by the definite positive matrices.
\end{itemize}

\section{Preliminaries}\label{prelim}

\subsection{Sublinear Operators}

We consider the space $\mathds R^N\times\mathds S^N$ with the inner product
\[
((p,S),(p',S')):=\frac12\langle S,S'\rangle+p^\dag p'
\]
and the norm $\|(p,S)\|:=\sqrt{((p,S),(p,S))}$.

\begin{assum}\label{subopassum}
In this subsection we will concentrate on the study of continuous operators of the form
\[
F:[0,T]\times\mathds R^N\times\mathds R^N\times\mathds S^N\to\mathds R
\]
such that the following properties hold true for any $(t,x)\in[0,T]\times\mathds R^N$, $(p,S)$ and $(p',S')$ in $\mathds R^N\times\mathds S^N$:
\begin{enumerate}[label=(\roman*)]
\item\label{subopassumsubad}\emph{(Subadditivity)} $F(t,x,p+p',S+S')\le F(t,x,p,S)+F(t,x,p',S')$;
\item\label{subopassumphom}\emph{(Positive Homogeneity)} If $\delta\ge0$ then $F(t,x,\delta p,\delta S)=\delta F(t,x,p,S)$;
\item\emph{(Uniform Ellipticity)} Exists a constant $\lambda>0$ such that, if $S'\ge0$,
\[
F(t,x,p,S+S')-F(t,x,p,S)\ge\lambda|S'|;
\]
\item\emph{(Lipschitz Continuity)} Exists a positive $\ell$ such that, for any $y\in\mathds R^N$,
\[
|F(t,x,p,S)-F(t,y,p,S)|\le\ell|x-y|\|(p,S)\|.
\]
\end{enumerate}
We will usually refer to $T$ as the terminal time of $F$, since it will play the role of terminal time in the parabolic problems which we will deal with later.
\end{assum}

The operators satisfying conditions \ref{subopassumsubad} and \ref{subopassumphom} are commonly known as \emph{sublinear operators}. Notice that \cref{subopassumphom,subopassumsubad} imply convexity in the third and fourth arguments and, vice versa, convexity and \ref{subopassumphom} imply \ref{subopassumsubad}.

The main result of this section is the following characterization theorem:

\begin{theo}\label{charsublin}
Let $F$ be as in \cref{subopassum} and $K_F$ be the set of the elements $(b,a)\in C^0\left([0,T]\times\mathds R^N;\mathds R^N\right)\times C^0\left([0,T]\times\mathds R^N;\mathds S_+^N\right)$ such that, for any $(t,x,p,S)\in[0,T]\times\mathds R^N\times\mathds R^N\times\mathds S^N$,
\[
\frac12\langle a(t,x),S\rangle+p^\dag b(t,x)\le F(t,x,p,S),
\]
$\Lip(b(t))\le2\ell$, $\Lip(a(t))\le2\sqrt2\ell$ and the eigenvalues of $a(t,x)$ are bigger than $2\lambda$. Then $K_F$ is a non empty and convex set, and
\[
F(t,x,p,S)=\max_{(b,a)\in K_F}\frac12\langle a(t,x),S\rangle+p^\dag b(t,x)
\]
for any $(t,x,p,S)\in[0,T]\times\mathds R^N\times\mathds R^N\times\mathds S^N$. Furthermore for each $(b,a)\in K_F$ the linear operator
\[
(t,x,p,S)\mapsto\frac12\langle a(t,x),S\rangle+p^\dag b(t,x)
\]
has the same ellipticity conditions of $F$.
\end{theo}

To prove this we preliminarily need the followings two lemmas. The first one is an adaptation of \cite[Lemma 1.8.14]{schneiderconvex}, which permit us to express the Hausdorff distance using support function, while the second one is just an adaptation of the Hahn--Banach theorem.

\begin{lem}\label{supfunhausdis}
Given, $A$ and $B$, two compact and convex subset of $\mathds R^N\times\mathds S^N$ we define the application
\[
h_A,\:h_B:\mathds R^N\times\mathds S^N\to\mathds R
\]
as the support functions of $A$ and $B$ respectively, that is to say
\[
h_A((p,S))=\sup_{(p',S')\in A}((p,S),(p',S')),\quad h_B((p,S))=\!\sup_{(p',S')\in B}((p,S),(p',S')),
\]
for any $(p,S)\in\mathds R^N\times\mathds S^N$. Then, for the Hausdorff distance
\begin{multline*}
d_H(A,B):=\\
\max\left\{\max_{(p,S)\in A}\min_{(p',S')\in B}\|(p,S)-(p',S')\|,\max_{(p',S')\in B}\min_{(p,S)\in A}\|(p,S)-(p',S')\|\right\},
\end{multline*}
we have that
\[
d_H(A,B)=\max_{\substack{(p,S)\in\mathds R^N\times\mathds S^N\\\|(p,S)\|=1}}|h_A((p,S))-h_B((p,S))|.
\]
\end{lem}

\begin{lem}\label{charsublinaux}
If $F$ is like in \cref{subopassum} then the set
\[
\mathcal L_{F(t,x)}:=\{L\in\left(\mathds R^N\times\mathds S^N\right)^*:L\le F(t,x)\}
\]
is non empty, compact and convex for any $(t,x)\in[0,T]\times\mathds R^N$. Furthermore $F(t,x)=\max\limits_{L\in\mathcal L_F(t,x)}L$.
\end{lem}

\proof[Proof of \cref{charsublin}]
By \cref{charsublinaux} and the Riesz representation theorem we have that, for any $(t,x)\in[0,T]\times\mathds R^N$, there exists the non empty convex and compact set
\[
K^x_t:=\left\{
\begin{gathered}
(b,a)\in\mathds R^N\times\mathds S^N:\frac12\langle a,S\rangle+p^\dag b\le F(t,x,p,S)\\
\text{for any }(p,S)\in\mathds R^N\times\mathds S^N
\end{gathered}
\right\}.
\]
Given a $(t,x)\in[0,T]\times\mathds R^N$ and $\left(\overline b,\overline a\right)\in K^x_t$ we define the function
\[
(b,a):[0,T]\times\mathds R^N\to\mathds R^N\times\mathds S^N
\]
such that
\[
(b,a)(s,y):=
\begin{cases}
\left(\overline b,\overline a\right),&\text{if }(s,y)=(t,x),\\
\argmin\limits_{(b,a)\in K_s^y}\left\|\left(\overline b-b,\overline a-a\right)\right\|,&\text{if }(s,y)\neq(t,x).
\end{cases}
\]
This function is well defined because is well known that the projection of a point onto a convex set, i.e. $\argmin\limits_{(b,a)\in K_s^y}\left\|\left(\overline b-b,\overline a-a\right)\right\|$, exists and is unique. We will show that $(b,a)\in K_F$, since this yields that $K_F$ is a non empty convex set (the convexity proof is trivial, hence we skip it) such that, thanks to the arbitrariness of the construction,
\[
F(t,x,p,S)=\max_{(b,a)\in K_F}\frac12\langle a(t,x),S\rangle+p^\dag b(t,x)
\]
for any $(t,x,p,S)\in[0,T]\times\mathds R^N\times\mathds R^N\times\mathds S^N$.\\
As a consequence of the definition and \cref{supfunhausdis} we have, for any $(s,y)$ in $[0,T]\times\mathds R^N$,
\begin{align*}
\|(b,a)(t,x)-(b,a)(s,y)\|=&\min_{(b,a)\in K_s^y}\|(b(t,x)-b,a(t,x)-a)\|\\
\le&\max_{(b_1,a_1)\in K_t^x}\min_{(b_2,a_2)\in K_s^y}\|(b_1-b_2,a_1-a_2)\|\\
\le&d_H(K^x_t,K^y_s)\\
=&\max_{\substack{(p,S)\in\mathds R^N\times\mathds S^N:\\\|(p,S)\|=1}}|F(t,x,p,S)-F(s,y,p,S)|.
\end{align*}
Since $|a|\le\sqrt2\|(b,a)\|$ for any $(b,a)\in\mathds R^N\times\mathds S^N$ and
\[
|a(r,z)-a(s,y)|\le|a(t,x)-a(s,y)|+|a(r,z)-a(t,x)|
\]
the previous inequality yields that $\Lip(a(s))\le2\sqrt2\ell$ for any $s\in[0,T]$, and similarly that $\Lip(b(s))\le2\ell$ for any $s\in[0,T]$.\\
We now prove the ellipticity part of the statement and, as a consequence, that $\sigma_{a(t,x)}\subset[2\lambda,\infty)^N$ for any $(t,x)\in[0,T]\times\mathds R^N$, thus that $(b,a)\in K_F$. Let $\lambda$ be the ellipticity constants of $F$ and
\begin{align*}
L:[0,T]\times\mathds R^N\times\mathds R^N\times\mathds S^N\to&\:\mathds R\\
(t,x,p,S)\mapsto&\:\frac12\langle a(t,x),S\rangle+p^\dag b(t,x),
\end{align*}
then, by its linearity, we only have to prove that for any $S\in\mathds S^N_+$ and $(t,x)$ in $[0,T]\times\mathds R^N$
\begin{equation}\label{eq:unelllineq}
L(t,x,0,S)\ge\lambda|S|.
\end{equation}
Obviously we have, for any $S\in\mathds S^N_+$ and $(t,x)\in[0,T]\times\mathds R^N$,
\[
\lambda|S|\le F(t,x,0,0)-F(t,x,0,-S)\le-L(t,x,0,-S)=L(t,x,0,S),
\]
hence \eqref{eq:unelllineq}. Finally, let $q$ an element of $\mathds R^N$ and define $Q:=qq^\dag$, which is an element of $\mathds S^N_+$ such that
\[
|Q|^2=\tr(qq^\dag qq^\dag)=\tr(q|q|^2q^\dag)=|q|^4.
\]
Therefore \eqref{eq:unelllineq} yields, for any $(t,x)\in[0,T]\times\mathds R^N$,
\[
\lambda|q|^2=\lambda|Q|\le\frac12\langle a(t,x),qq^\dag\rangle=\frac12q^\dag a(t,x)q
\]
and the Rayleigh quotient formula proves that $\sigma_{a(t,x)}\subset[2\lambda,\infty)^N$, concluding the proof.
\endproof

\begin{rem}
We point out that can be easily proved that for the previous theorem holds a converse.
\end{rem}

We have characterized $F$ as the support function of a set of linear operators. Usually, to obtain representation formulas for viscosity solutions to a second order PDE with linear operator like \eqref{exparsemlin}, is useful to study a function $\sigma$ such that $\sigma\sigma^\dag$ is the diffusion part of that operator, hence we will do something similar: if we define the application from $\mathds S^N_+$ to itself which associate via singular value decomposition the matrix $a$ with its square root $\sigma$ then it is well defined, as can be seen in \cite[Section 6.5]{bellmanmatrixbook}. Moreover we know from \cite[Lemma 2.1]{sqmatrix} that, on the space of matrices with eigenvalues equal or bigger than $2\lambda$, this application is Lipschitz continuous with Lipschitz constant $c_\lambda:=\frac1{2\sqrt\lambda}$, therefore the application $(b,\sigma)\mapsto(b,\sigma\sigma)$ that maps the set $\mathcal K_F$, which contains the $(b,\sigma)$ in $C^0\left([0,T]\times\mathds R^N;\mathds R^N\right)\times C^0\left([0,T]\times\mathds R^N;\mathds S^N_+\right)$ such that $(b,\sigma\sigma)\in K_F$ and $\Lip(\sigma(t))\le2\sqrt2c_\lambda\ell$ for any $t\in[0,T]$,
into $K_F$ is surjective and consequently
\[
F(t,x,p,S)=\max_{(b,\sigma)\in\mathcal K_F}\frac12\left\langle\sigma^2(t,x),S\right\rangle+p^\dag b(t,x).
\]

Our method to obtain representation formulas relies on a \emph{dynamic programming principle}, which will be presented later and is based on a construction on a broader set than $\mathcal K_F$. This set, which we call $\mathcal A_F$, is made up of the functions
\[
(b,\sigma):[0,T]\times\Omega\times\mathds R^N\to\mathds R^N\times\mathds R^{N\times N}
\]
which are cadlag, i.e. right continuous and left bounded, on $[0,T]$ and such that, for any $(t,x,p,S)\in[0,T]\times\mathds R^N\times\mathds R^N\times\mathds S^N$ and $\omega\in\Omega$,
\[
\frac12\left\langle\sigma^2(t,\omega,x),S\right\rangle+p^\dag b(t,\omega,x)\le F(t,x,p,S),
\]
$\Lip(b(t,\omega))\le2\ell$, $\Lip(\sigma(t,\omega))\le2\sqrt 2c_\lambda\ell$, the eigenvalues of $(\sigma\sigma)(t,\omega,x)$ belong to $[2\lambda,\infty)^N$ and $\{(b,\sigma)(t,0)\}_{t\in[0,T]}$ is a progressive process. $\mathcal A_F$ is obviously non empty, since it contains $\mathcal K_F$. For any stopping time $\tau$ with value in $[0,T]$, an useful subset of $\mathcal A_F$, which we will use later, is $\mathcal A^\tau_F$, which consists of the $(b,\sigma)$ belonging to $\mathcal A_F$ such that $\{(b,\sigma)(\tau+t,x)\}_{t\in[0,\infty)\!}$ is progressive with respect to the filtration $\!\left\{\mathcal F_t^\tau\right\}_{t\in[0,\infty)}$. Trivially $\mathcal A^0_F=\mathcal A_F$.

\begin{rem}\label{aflp}
It is easy to see that, for any $p>0$, the process $(b,\sigma)(s,0)$ belongs to $\mathds L_N^p(t)\times\mathds L^p_{N\times N}(t)$ for each $(b,\sigma)\in\mathcal A_F$ and $t\in[0,\infty)$, since the image of $[0,t]\times\Omega\times\{0\}$ under $(b,\sigma)$ is contained on a compact set for any $(b,\sigma)\in\mathcal A_F$ and $t\in[0,\infty)$.
\end{rem}

\subsection{Dynamic Programming Principle}\label{dynprincsec}

The scope of this section is to provide the necessary instruments to prove the \cref{cauchydpp}, which will be used to derive representation formulas for viscosity solutions to the parabolic problem which we will study later. The dynamic programming principle is, in this contest, an instrument that permit us to break a stochastic trajectory in two or more part. In particular for the problem \eqref{exparsemlin} it means that
\begin{multline*}
\mathds E\left(u\left(s,X^{t,x}_s\right)+\int_s^Tf\left(r,X_r^{s,X^{t,x}_s},Y_r^{s,X^{t,x}_s},Z_r^{s,X^{t,x}_s}\right)dr\right)\\
\begin{aligned}
=&\mathds E\left(Y^{s,X^{t,x}_s}_s+\int_s^Tf\left(r,X_r^{s,X^{t,x}_s},Y_r^{s,X^{t,x}_s},Z_r^{s,X^{t,x}_s}\right)dr\right)\\
=&\mathds E\left(Y^{t,x}_s+\int_s^Tf\left(r,X_r^{t,x},Y_r^{t,x},Z_r^{t,x}\right)dr\right)\\
=&\mathds E\left(Y^{t,x}_t\right)=u(t,x),
\end{aligned}
\end{multline*}
which is just a simple consequence of the uniqueness of the solutions to the FBSDE, while for the $G$--heat equation \eqref{gheat} this means that
\begin{multline}\label{gheatdpp}
\sup_{\sigma\in\mathcal A_G}\mathds E\left(u\left(s,x+\int_t^s\sigma(r)dW_r\right)\right)\\
\begin{aligned}
=&\sup_{\sigma\in\mathcal A_G}\sup_{\sigma'\in\mathcal A_G}\mathds E\left(g\left(x+\int_t^s\sigma(r)dW_r+\int_s^T\sigma'(r)dW_r\right)\right)\\
=&\sup_{\sigma\in\mathcal A_G}\mathds E\left(g\left(x+\int_t^T\sigma(r)dW_r\right)\right)\\
=&u(t,x).
\end{aligned}
\end{multline}
The proof of \eqref{gheatdpp} is contained in \cite[Subsection 3.1]{art:pathpeng}. This also intuitively explain why we ask to the elements of $\mathcal A_F$ to be cadlag in time. We point out that in \cite{art:pathpeng} the authors ask to the elements of $\mathcal A_F$ to only be measurable in time, but for the analysis of the more general problem \eqref{parsubopgen} in \cref{parprob} we will use right continuity.

The dynamic programming principle exposed in theorem \ref{cauchydpp} is a generalization of the one presented by Denis, Hu and Peng in \cite{art:pathpeng} and can be obtained using a similar method, slightly adapting the proofs of \cite[Lemmas 41--44]{art:pathpeng}. Hence we will just present the results that are relevant for our analysis skipping the proofs.

First of all we notice that for any $(b,\sigma)\in\mathcal A_F$ we can define an $X_{(b,\sigma)}$ solution to the SDE $(b,\sigma)$ as in \eqref{eq:sde}. However, to ease notations, we will usually omit the dependence of $X$ from $(b,\sigma)$.

Given a sublinear operator $F$ with a positive terminal time $T$ as in \cref{subopassum}, a stopping time $\tau$ with value in $[0,T]$ and a measurable application $\varphi$ from $[0,T]\times\mathds R^N\times\mathcal A_F$ into $\mathds R$, continuous in probability with respect to $\mathcal A_F$ and such that
\[
\varphi_\tau:\mathds R^N\times\mathcal A_F|_{[\tau,T]}\to L^1\left(\Omega,\mathcal F_T;\mathds R\right),
\]
we define, for any $\zeta\in L^2\left(\Omega,\mathcal F_\tau;\mathds R^N\right)$, the function
\[
\Phi_\tau(\zeta):=\esss_{(b,\sigma)\in\mathcal A_F}\mathds E(\varphi_\tau(\zeta,b,\sigma)|\mathcal F_\tau).
\]
We assume that $\varphi_t(x,b,\sigma)$ is $\mathcal F^t_\infty$--measurable for any $t\in[0,T]$, $x\in\mathds R^N$ and $(b,\sigma)\in\mathcal A_F^t$, and, for any stopping time $\tau$ with value in $[0,T]$ and $\zeta$ in $L^2\left(\Omega,\mathcal F_\tau;\mathds R^N\right)$,
\begin{equation}\label{conddpp}
\sup_{(b,\sigma)\in\mathcal A_F}\mathds E(|\varphi_\tau(\zeta,b,\sigma)|)<\infty.
\end{equation}

In this section the function $\Phi$ represents, roughly speaking, the viscosity solution, $\zeta$ is the first part of a stochastic trajectory broken off at $\tau$ (this is why we restrict $\varphi_\tau$ on $\mathcal A_F|_{[\tau,T]}$) and $\varphi$ the function which we will use to build the viscosity solution.

\begin{lem}\label{peng40}
For each $(b_1,\sigma_1),(b_2,\sigma_2)$ in $\mathcal A_F$ and stopping time $\tau$ with value in $[0,T]$ there exists an $(b,\sigma)\in\mathcal A_F$ such that
\[
\mathds E(\varphi_\tau(\zeta,b,\sigma)|\mathcal F_\tau)=\mathds E(\varphi_\tau(\zeta,b_1,\sigma_1)|\mathcal F_\tau)\vee\mathds E(\varphi_\tau(\zeta,b_2,\sigma_2)|\mathcal F_\tau).
\]
Therefore there exists a sequence $\{(b_i,\sigma_i)\}_{i\in\mathds N}$ in $\mathcal A_F$ such that a.e.
\[
\mathds E(\varphi_\tau(\zeta,b_i,\sigma_i)|\mathcal F_\tau)\uparrow\Phi_\tau(\zeta).
\]
We also have
\[
\mathds E(|\Phi_\tau(\zeta)|)\le\sup_{(b,\sigma)\in\mathcal A_F}\mathds E(|\varphi_\tau(\zeta,b,\sigma)|)<\infty,
\]
and, for any stopping time $\tau'\le\tau$,
\[
\mathds E\left(\esss_{(b,\sigma)\in\mathcal A_F}\mathds E(\varphi_\tau(\zeta,b,\sigma)|\mathcal F_\tau)\middle|\mathcal F_{\tau'}\right)=\esss_{(b,\sigma)\in\mathcal A_F}\mathds E(\varphi_\tau(\zeta,b,\sigma)|\mathcal F_{\tau'}).
\]
\end{lem}

\begin{rem}
To prove \cref{peng40} the randomness of the elements of $\mathcal A_F$ is crucial, this is the reason why we consider a set of stochastic process instead of a deterministic one.
\end{rem}

To continue we need a density result on $\mathcal A_F$ endowed with the topology of the $L^2$--convergence on compact set, which is to say that a sequence in $\mathcal A_F$ converges to an element of $\mathcal A_F$ if and only if it converges in $L^2([0,T]\times\Omega\times K)$ for any compact set $K\subset\mathds R^N$.

\begin{lem}
The set
\[
\mathcal J^t:=\left\{
\begin{aligned}
\alpha\in\mathcal A_F:\alpha|_{[t,T]}=\sum_{i=0}^n\chi_{A_i}\alpha_i|_{[t,T]},\text{ where }\{\alpha_i\}_{i=0}^n\subset\mathcal A_F^t\\
\text{and }\{A_i\}_{i=0}^n\text{ is a }\mathcal F_t\text{--partition of }\Omega
\end{aligned}
\right\}
\]
is dense in $\mathcal A_F$ for any $t\in[0,T]$.
\end{lem}
\proof
To prove this we will show that, fixed a $k\in\mathds N$, we can approximate, in $L^2([0,T]\times\Omega\times B_k(0))$, any element of $\mathcal A_F$ with an element of $\mathcal J^t$.\\
Preliminarily notice that by our assumption each element of $\mathcal A_F$ can be approximated in $L^2([0,T]\times\Omega\times B_k(0))$ by a sequence of simple functions. We will denote with $\mathcal B([0,T]\times B_k(0))$ the Borel $\sigma$--algebra of $[0,T]\times B_k(0)$.\\
Furthermore, since the collection $\mathcal I$ of the rectangles $A\times B$ where $A\in\mathcal F_T$ and $B\in\mathcal B([0,T]\times B_k(0))$ is a $\pi$--system which contains the complementary of its sets and generate $\sigma(\mathcal F_T\times\mathcal B([0,T]\times B_k(0)))$, by \mbox{\cite[Dynkin's lemma A1.3]{williams08}} each set in $\sigma(\mathcal F_T\times\mathcal B([0,T]\times B_k(0)))$, which is the smallest $d$--system containing $\mathcal I$, can be approximate by a finite union of sets in $\mathcal I$. Similarly, each set in $\mathcal F_T$ can be approximated by finite intersection and union of sets in $\mathcal F_t$ and $\mathcal F_T^t$, since $\mathcal F_T=\sigma\left(\mathcal F_t,\mathcal F_T^t\right)$.\\
Therefore, fixed $(b,\sigma)\in\mathcal A_F$, for any $\varepsilon>0$ there exists a simple function $s_\varepsilon$ such that $s_\varepsilon(t,\omega,x)=\sum\limits_{i=1}^n\sum\limits_{j=1}^ms^j_i(t,x)\chi_{A_i}(\omega)\chi_{A'_j}(\omega)$ where $\{A_i\}_{i=1}^n$ and $\{A'_j\}_{j=1}^m$ are respectively a $\mathcal F_t$--partition and a $\mathcal F^t_T$--partition of $\Omega$ and
\begin{equation}\label{eq:exJdenseA1}
\mathds E\left(\int_0^T\int_{B_k(0)}|(b,\sigma)(t,x)-s_\varepsilon(t,x)|^2dxdt\right)<\varepsilon.
\end{equation}
Then, for each $A_i$ and $A'_j$ with $\mathds P(A_i\cap A'_j)>0$, there exists a $\omega_i^j\in A_i\cap A'_j$ such that
\[
\int_0^T\int_{B_k(0)}\left|(b,\sigma)\left(t,\omega_i^j,x\right)-s_i^j(t,x)\right|^2dxdt<\frac\varepsilon{\mathds P(A_i\cap A'_j)},
\]
otherwise we would have that
\begin{multline*}
\mathds E\left(\int_0^T\int_{B_k(0)}|(b,\sigma)(t,x)-s_\varepsilon(t,x)|^2dxdt\right)\\
\begin{aligned}
=&\mathds E\left(\sum_{i=1}^n\sum_{j=1}^m\int_0^T\int_{B_k(0)}\left|(b,\sigma)(t,x)-s_i^j(t,x)\right|^2\chi_{A_i}\chi_{A'_j}dtdx\right)\\
\ge&\mathds E\left(\int_0^T\int_{B_k(0)}\left|(b,\sigma)(t,x)-s_i^j(t,x)\right|^2\chi_{A_i}\chi_{A'_j}dtdx\right)\\
\ge&\varepsilon,
\end{aligned}
\end{multline*}
in contradiction with \eqref{eq:exJdenseA1}. Finally, let $\omega_i^j$ be any elements of $A_i\cap A'_j$ if $\mathds P(A_i\cap A'_j)=0$ and $\left(b_i^k,\sigma_i^k\right):=\sum\limits_{j=1}^m(b,\sigma)\left(\omega_i^j\right)\chi_{A'_j}$. Then $\left(b_i^k,\sigma_i^k\right)\in\mathcal A^t_F$ and $\left(b^k_\varepsilon,\sigma^k_\varepsilon\right):=\sum\limits_{i=1}^n\left(b^k_i,\sigma^k_i\right)\chi_{A_i}$ is an element of $\mathcal J^t$ satisfying
\[
\mathds E\left(\int_0^T\int_{B_k(0)}\left|(b,\sigma)(t,x)-\left(b^k_\varepsilon,\sigma^k_\varepsilon\right)(t,x)\right|^2dxdt\right)<4\varepsilon.
\]
This proves that $\mathcal J^t$ is dense in $\mathcal A_F$.
\endproof

\begin{lem}
For each $t\in[0,T]$ and $x\in\mathds R^N$, $\Phi_t(x)$ is deterministic. Furthermore
\begin{equation}\label{eq:artpeng42}
\Phi_t(x)=\esss_{(b,\sigma)\in\mathcal A_F}\mathds E(\varphi_t(x,b,\sigma)|\mathcal F_t)=\esss_{(b,\sigma)\in\mathcal A_F^t}\mathds E(\varphi_t(x,b,\sigma)|\mathcal F_t).
\end{equation}
\end{lem}

\begin{lem}\label{artpeng43}
We define the function
\begin{align*}
u:[0,T]\times\mathds R^N&\xrightarrow{\qquad}\mathds R\\
(t,x)&\xmapsto{\qquad}\Phi_t(x)
\end{align*}
and assume that it is continuous. Then, for each stopping time $\tau$ with value in $[0,T]$ and $\zeta\in L^2_N\left(\Omega,\mathcal F_\tau;\mathds R^N\right)$, we have that $u_\tau(\zeta)=\Phi_\tau(\zeta)$ a.e..
\end{lem}

\begin{rem}\label{supeqrem}
This lemma says, as a consequence of \eqref{eq:artpeng42}, that
\[
\esss_{(b,\sigma)\in\mathcal A_F}\mathds E(\varphi_t(\zeta,b,\sigma)|\mathcal F_t)=\esss_{(b,\sigma)\in\mathcal A_F^t}\mathds E(\varphi_t(\zeta,b,\sigma)|\mathcal F_t)
\]
for any $t\in[0,T]$.
\end{rem}

\section{Parabolic PDEs with Sublinear Operators}\label{parprob}

We analyze now the following problem:
\begin{prob}\label{parsubopgenprob}
Let $T$ be a terminal time, $F$ a uniformly elliptic operator satisfying \cref{subopassum} and
\[
f:[0,T]\times\mathds R^N\times\mathds R\times\mathds R^N\to\mathds R\text{ and }g:\mathds R^N\to\mathds R
\]
two continuous functions for which there exist two constants $\mu\in\mathds R$ and $\ell\ge0$ such that, for any $t\in[0,T]$, $x,x'\in\mathds R^N$, $y,y'\in\mathds R$ and $z,z'\in\mathds R^N$,
\begin{enumerate}[label=(\roman*)]
\item $|g(x)-g(x')|\le\ell|x-x'|$;
\item $|g(x)|\le\ell(1+|x|)$;
\item $|f(t,x,y,z)-f(t,x',y,z')|\le\ell(|x-x'|+|z-z'|)$;
\item $|f(t,x,y,z)|\le\ell(1+|x|+|y|+|z|)$;
\item $(y-y')(f(t,x,y,z)-f(t,x,y',z))\le\mu|y-y'|^2$.
\end{enumerate}
Find the solution $u$ to the parabolic PDE
\[
\begin{cases}
\partial_tu(t,x)+F\left(t,x,\nabla_xu,D^2_xu\right)+f(t,x,u,\nabla_xu)=0,&t\in(0,T),x\in\mathds R^N,\\
u(T,x)=g(x),&x\in\mathds R^N.
\end{cases}
\]
\end{prob}

\begin{rem}
To ease notations, we can assume without loss of generality that the $\ell$ in \cref{parsubopgenprob} is the same as in \cref{subopassum}. Since $F$ is continuous, we can also assume that, for any $(b,\sigma)\in\mathcal A_F$, $|b(t,0)|\le\ell$ and $|\sigma(t,0)|\le\ell$ for any $t\in[0,T]$. We also define the sets $\mathcal L_F^\tau$ whose elements are the operators $L_{(b,\sigma)}$ with $(b,\sigma)\in\mathcal A^\tau_F$ such that
\[
L_{(b,\sigma)}(p,S):=\frac12\left\langle\sigma^2(t,x),S\right\rangle+p^\dag b(t,x).
\]
As previously done with $\mathcal A_F$, we also define the set $\mathcal L_F:=\mathcal L_F^0$.
\end{rem}

Let us define what we mean with viscosity solution to \cref{parsubopgenprob}. For a detailed overview of the viscosity solution theory we refer to \cite{userguide}.

\begin{defin}
Given an upper semicontinuous function $u$ we say that a function $\varphi$ is a \emph{supertangent} to $u$ at $(t,x)$ if $(t,x)$ is a local maximizer of $u-\varphi$.\\
Similarly we say that a function $\psi$ is a \emph{subtangent} to a lower semicontinuous function $v$ at $(t,x)$ if $(t,x)$ is a local minimizer of $v-\psi$.
\end{defin}

\begin{defin}
An upper semicontinuous function $u$ is called a \emph{viscosity subsolution} to \cref{parsubopgenprob} if, for any suitable $(t,x)$ and $C^{1,2}$ supertangent $\varphi$ to $u$ at $(t,x)$,
\[
\partial_t\varphi(t,x)+F\left(t,x,\nabla_x\varphi,D^2\varphi\right)+f(t,x,u,\nabla_x\varphi)\ge0.
\]
Similarly a lower semicontinuous function $v$ is called a \emph{viscosity supersolution} to \cref{parsubopgenprob} if, for any suitable $(t,x)$ and $C^{1,2}$ subtangent $\psi$ to $v$ at $(t,x)$,
\[
\partial_t\psi(t,x)+F\left(t,x,\nabla_x\psi,D^2\psi\right)+f(t,x,u,\nabla_x\psi)\le0.
\]
Finally a continuous function $u$ is called a \emph{viscosity solution} to \cref{parsubopgenprob} if it is both a super and a subsolution to \cref{parsubopgenprob}.
\end{defin}

We derive from \cref{defviscomp} that a comparison result holds true.

\begin{theo}\label{parsubopgencomp}
Let $u$ and $v$ be respectively a subsolution and a supersolution to \cref{parsubopgenprob} satisfying polynomial growth condition. If $u|_{t=T}\le v|_{t=T}$, then $u\le v$ on $(0,T]\times\mathds R^N$.
\end{theo}

When $F$ is a linear operator it is known that the representation formula of its viscosity solution is built from a FBSDE, see \eqref{exparsemlin}. To adapt this method to our case we will use the \cref{cauchydpp}.

\begin{defin}\label{fbsde}
Consider the FBSDE
\begin{equation}\label{eq:fbsde}
\left\{\begin{aligned}
X^{t,\zeta}_s=&\zeta+\int_t^s\sigma\left(r,X^{t,\zeta}_r\right)dW_r+\int_t^sb\left(r,X^{t,\zeta}_r\right)dr,\\
Y^{t,\zeta}_s=&g\left(X^{t,\zeta}_T\right)+\int_s^Tf_\sigma\left(r,X_r^{t,\zeta},Y_r^{t,\zeta},Z_r^{t,\zeta}\right)dr\\
&-\int_s^TZ_r^{t,\zeta}dW_r,
\end{aligned}\right.\,s\in[t,T],
\end{equation}
where $\zeta\in L^2\left(\Omega,\mathcal F_t;\mathds R^N\right)$, $(b,\sigma)\in\mathcal A_F$, the function $f_\sigma$ is defined as
\[
f_\sigma(t,x,y,z):=f\left(t,x,y,z(\sigma(t,x))^{-1}\right),
\]
for any $(t,x,y,z)$ in $[0,T]\times\mathds R^N\times\mathds R\times\mathds R^N$ and the functions $f$ and $g$ are as in the assumptions of \cref{parsubopgenprob}. Thanks to the uniformly ellipticity condition we know that $f_\sigma$ is well defined and that the Lipschitz constant for its the fourth argument is $\ell\sqrt{\frac N{2\lambda}}$, but for simplicity we will assume that it is $\ell$, possibly increasing it.\\
Note that under this conditions the \cref{sdeassum,bsdeassum} hold for $X$ and $(Y,Z)$ respectively. We will call $(X,Y,Z)$ a solution to the FBSDE if $X$ is a solution to the SDE part of this system and $(Y^{t,\zeta},Z^{t,\zeta})$ is a solution to the BSDE part for any $(t,\zeta)\in[0,T]\times L^2\left(\Omega,\mathcal F_t;\mathds R^N\right)$. Notice that, under our assumptions, there exists a unique solution to \eqref{eq:fbsde}, thanks to \cref{sdeex,bsdeex}. Due to \cref{sdestopstart,stopbsde}, this is true even if $t$ is an a.e. finite stopping time.
\end{defin}

\begin{rem}
Notice that the uniqueness property of the FBSDE imply that, for any $0\le t\le r\le s\le T$,
\[
\left(X^{r,X^{t,\zeta}_r}_s,Y^{r,X^{t,\zeta}_r}_s,Z^{r,X^{t,\zeta}_r}_s\right)=\left(X^{t,\zeta}_s,Y^{t,\zeta}_s,Z^{t,\zeta}_s\right).
\]
This holds true even if $t$, $r$ and $s$ are stopping time.
\end{rem}

\begin{rem}\label{L_FtoI_F}
We point out that since the elements of $\mathcal L_F^\tau$ and the solutions to SDEs $(b,\sigma)$ can be uniquely determined, except for the initial data of the SDEs, by an element of $\mathcal A_F^\tau$, we can uniquely link to each operator $L\in\mathcal L_F^\tau$ an $X_{(b,\sigma)}$. Moreover, for each problem, we can uniquely associate in the same way a solution of the FBSDE \eqref{eq:fbsde}.
\end{rem}

For the remainder of this section, we will simply write $Y$ to denote the second term of the triplet $(X,Y,Z)$ solution to the FBSDE defined in \cref{fbsde}, for $(b,\sigma)$ that varies in $\mathcal A_F$. For simplicity we will omit the dependence of $Y$ and $Z$ from $X$ and $\sigma$ or, equivalently, from $(b,\sigma)$.

We will prove that $u(t,x):=\sup\limits_{(b,\sigma)\in\mathcal A_F}\mathds E\left(Y^{t,x}_t\right)$ is a viscosity solution to the \cref{parsubopgenprob} breaking the proof in several steps.

\begin{prop}\label{parsubopgenviscont}
The function $u(t,x):=\sup\limits_{(b,\sigma)\in\mathcal A_F}\mathds E\left(Y^{t,x}_t\right)$ is $\frac12$--Hölder continuous in the first variable and Lipschitz continuous in the second one. Furthermore we have that there exists a constant $c$, which depends only on $\ell$, $\mu$ and $T$, such that
\begin{equation}\label{eq:parvisbound}
\mathds E\left(|u(\tau,\zeta)|^2\right)\le\sup_{(b,\sigma)\in\mathcal A_F}\mathds E\left(\left|Y^{\tau,\zeta}_\tau\right|^2\right)\le c\left(1+\mathds E\left(|\zeta|^2\right)\right),
\end{equation}
for any $t\in[0,T]$ and $\zeta\in L^2\left(\Omega,\mathcal F_t;\mathds R^N\right)$.
\end{prop}

We point out that this proposition permits us to use the results of \cref{dynprincsec} on $u$. In particular $Y^\tau_\tau$, which is $\mathcal F^\tau_\tau$--measurable and therefore a.e. deterministic, has the same role of $\varphi$ in \cref{dynprincsec}. We already know that $Y^\tau_\tau$ is continuous in probability, thanks to our assumptions, \cref{sdeex,bsdeex}, furthermore we prove here that it satisfies \eqref{conddpp} and the continuity of $u$, which is needed for \cref{artpeng43}.

\proof
To prove our statement note that by the definition and the Jensen's inequality
\begin{align*}
|u(t,x)-u(s,y)|=&\left|\sup_{(b,\sigma)\in\mathcal A_F}\mathds E\left(Y^{t,x}_t\right)-\sup_{(b,\sigma)\in\mathcal A_F}\mathds E\left(Y^{s,y}_s\right)\right|\\
\le&\sup_{(b,\sigma)\in\mathcal A_F}\left(\mathds E\left(\left|Y^{t,x}_t-Y^{s,y}_s\right|^2\right)\right)^\frac12
\end{align*}
and
\[
|u(t,\zeta)|=\left|\sup_{(b,\sigma)\in\mathcal A_F}\mathds E\left(Y^{t,\zeta}_t\right)\right|\le\left(\sup_{(b,\sigma)\in\mathcal A_F}\mathds E\left(\left|Y^{t,\zeta}_t\right|^2\right)\right)^\frac12,
\]
for any $t,s\in[0,T]$, $x,y\in\mathds R^N$ and $\zeta\in L^2\left(\Omega,\mathcal F_t;\mathds R^N\right)$. The statement is then a consequence of \cref{bsdeex,sdeex}.
\endproof

We can now prove the dynamic programming principle for $u$.

\begin{theo}[Dynamic programming principle]\label[dpp]{cauchydpp}
For any $(b,\sigma)\in\mathcal A_F$ we let $\left(\overline Y,\overline Z\right)$ be the solution of the BSDE
\begin{equation}\label{eq:cauchydpp1}
\overline Y_s=u\left(\tau,X_\tau^{t,x}\right)+\int_{s\wedge\tau}^\tau f_\sigma\left(r,X^{t,x}_r,\overline Y_r,\overline Z_r\right)dr-\int_{s\wedge\tau}^\tau\overline Z_rdW_r,
\end{equation}
where $s\in[t,T]$ and $\tau$ is a stopping time with value in $[t,T]$. Then we have $\sup\limits_{(b,\sigma)\in\mathcal A_F}\mathds E\left(\overline Y_t\right)=u(t,x)$.
\end{theo}
\proof
Fix $(\overline b,\overline\sigma)\in\mathcal A_F$ in \eqref{eq:cauchydpp1} and define $\overline X:=X_{(b,\sigma)}$ and the subset of $\mathcal A_F$
\[
\overline{\mathcal A}_F:=\left\{(b,\sigma)\in\mathcal A_F:(b,\sigma)(s)=\left(\overline b,\overline\sigma\right)(s)\text{ for any }s\in[t,\tau]\right\}.
\]
From \cref{artpeng43} we know that
\[
\esss_{(b,\sigma)\in\overline{\mathcal A}_F}\mathds E\left(Y_\tau^{t,x}\middle|\mathcal F_\tau\right)=\esss_{(b,\sigma)\in\overline{\mathcal A}_F}\mathds E\left(Y_\tau^{\tau,\overline X^{t,x}_\tau}\middle|\mathcal F_\tau\right)=u\left(\tau,\overline X_\tau^{t,x}\right)
\]
and \cref{peng40} yields the existence of a sequence $\{(b_n,\sigma_n)\}_{n\in\mathds N}$ in $\overline{\mathcal A}_F$ and a corresponding sequence $\{Y_n\}_{n\in\mathds N}$ such that
\[
\lim_{n\to\infty}\mathds E\left(Y_{n,\tau}^{t,x}\middle|\mathcal F_\tau\right)=\esss_{(b,\sigma)\in\overline{\mathcal A}_F}\mathds E\left(Y_\tau^{t,x}\middle|\mathcal F_\tau\right)=u\left(\tau,\overline X_\tau^{t,x}\right).
\]
Then, by \cref{bsdeex} and the dominated convergence theorem, there exists a constant $c$ such that
\[
\lim_{n\to\infty}\mathds E\left(\left|\overline Y_t-Y^{t,x}_{n,t}\right|^2\right)\le\lim_{n\to\infty}c\mathds E\left(\left|u\left(\tau,\overline X_\tau^{t,x}\right)-Y^{t,x}_{n,\tau}\right|^2\right)=0,
\]
hence, up to subsequences,
\begin{equation}\label{eq:cauchydpp2}
\lim\limits_{n\to\infty}\mathds E\left(Y_{n,t}^{t,x}\right)=\mathds E\left(\overline Y_t\right).
\end{equation}
Furthermore, thanks to \cref{bsdeconf}, $Y^{t,x}_t\le\overline Y_t$ for any $(b,\sigma)\in\overline{\mathcal A}_F$, which together with \eqref{eq:cauchydpp2} implies that $\sup\limits_{(b,\sigma)\in\overline{\mathcal A}_F}\mathds E\left(Y^{t,x}_t\right)=\mathds E\left(\overline Y_t\right)$. Therefore we can use the arbitrariness of $\left(\overline b,\overline\sigma\right)$ to obtain our conclusion:
\[
\sup_{(b,\sigma)\in\mathcal A_F}\mathds E\left(\overline Y_t\right)=\sup_{(b,\sigma)\in\mathcal A_F}\mathds E\left(Y_t^{t,x}\right)=u(t,x).
\]
\endproof

Now we proceed to show that $u$ is a viscosity subsolution. In order to do that, we need the following lemma:

\begin{lem}\label{parsemlinvis}
For any $t\in(0,T)$, let $L$ be an element of $\mathcal L^t_F$ and $(X,Y,Z)$ the solution to the FBSDE \eqref{eq:fbsde} associated to $L$ as in \cref{L_FtoI_F}. If we define, for any $x\in\mathds R^N$ and $s\in[t,T]$, $u_L(s,x):=\mathds E\left(Y^{s,x}_s\right)$ we have that, for any supertangent $\varphi$ to $u_L$ at $(t,x)$,
\[
L\left(t,x,\nabla_x\varphi,D^2_x\varphi\right)\ge-\partial_t\varphi(t,x)-f\left(t,x,u_L,\nabla_x\varphi\right).
\]
\end{lem}
\proof
We preliminarily denote by $(b,\sigma)$ the element of $\mathcal A_F^t$ associated to $L$ and point out that since $(b,\sigma)$, restricted in $[t,T]$, is progressive with respect to the $\sigma$--algebra $\{\mathcal F^t_s\}_{s\in[t,T]}$, so are $L$, $X^t$ and $Y^t$. They are therefore constants a.e. in $t$. As a consequence $u_L(t,x)=Y^{t,x}_t$ a.e. for any $x\in\mathds R^N$.\
Given $x\in\mathds R^N$ and a supertangent $\varphi$ to $u_L$ at $(t,x)$ we can assume without loss of generality that $u_L(t,x)=\varphi(t,x)$, so we suppose that, a.e.,
\begin{equation}\label{eq:semlinvis0}
\partial_t\varphi(t,x)+L\left(t,x,\nabla_x\varphi,D^2_x\varphi\right)+f_\sigma(t,x,u_L,\nabla_x\varphi\sigma)<0
\end{equation}
and we will find a contradiction. Note that, as a consequence of the Blumenthal's 0--1 law, \eqref{eq:semlinvis0} is a deterministic inequality a.e.. By the definition of supertangent, there exists a $\delta\in(0,T-t)$ such that, for any $s\in[t,t+\delta]$ and $y\in B_\delta(x)$,
\begin{equation}\label{eq:semlinvis1}
u_L(s,y)\le\varphi(s,y),
\end{equation}
hence we define the stopping time
\[
\tau:=(t+\delta)\wedge\inf\left\{s\in[t,\infty):\left|X_s^{t,x}-x\right|\ge\delta\right\}
\]
and assume, possibly taking a smaller $\tau$, that
\begin{equation}\label{eq:semlinvis2}
\begin{aligned}
\partial_t\varphi\left(s\wedge\tau,X^{t,x}_{s\wedge\tau}\right)+L\left(s\wedge\tau,X^{t,x}_{s\wedge\tau},\nabla_x\varphi,D^2_x\varphi\right)&\\
+f_\sigma\left(s\wedge\tau,X^{t,x}_{s\wedge\tau},\varphi,\nabla_x\varphi\sigma\right)&<0.
\end{aligned}
\end{equation}
We point out that, by \eqref{eq:semlinvis0} and \cref{cadlagcont}, the previous inequality holds true on a set of positive measure for the $\chi_{[t,t+\delta)}dt\times d\mathds P$ measure, thus $\tau>t$ on a set of positive measure.\\
Let $\left(\overline Y_s,\overline Z_s\right):=\left(Y^{t,x}_{s\wedge\tau},Z^{t,x}_{s\wedge\tau}\right)$, which solve the BSDE
\[
\overline Y_s=Y_\tau^{t,x}+\int_{s\wedge\tau}^\tau f_\sigma\left(r,X^{t,x}_r,\overline Y_r,\overline Z_r\right)dr-\int_{s\wedge\tau}^\tau\overline Z_rdW_r,\quad s\in[t,T],
\]
and $\left(\hat Y_s,\hat Z_s\right):=\left(\varphi\left(s,X^{t,x}_{s\wedge\tau}\right),(\nabla_x\varphi\sigma)\left(s,X^{t,x}_{s\wedge\tau}\right)\right)$ which, by Itô's formula, is solution to
\[
\begin{aligned}
\hat Y_s=&\varphi\left(\tau,X^{t,x}_\tau\right)-\int_{s\wedge\tau}^\tau\hat Z_rdW_r\\
&-\int_{s\wedge\tau}^\tau\left(\partial_t\varphi\left(r,X^{t,x}_r\right)+L\left(r,X^{t,x}_r,\nabla_x\varphi,D^2_x\varphi\right)\right)dr,
\end{aligned}
\quad s\in[t,T].
\]
By \eqref{eq:semlinvis1} we have that
\[
u_L\left(\tau,X_\tau^{t,x}\right)-\varphi\left(\tau,X_\tau^{t,x}\right)=Y_\tau^{\tau,X_\tau^{t,x}}-\varphi\left(\tau,X_\tau^{t,x}\right)\le0
\]
and \eqref{eq:semlinvis2} imply, thanks to \cref{bsdeconf}, that $Y_t^{t,x}<\varphi(t,x)$ a.e., but this lead to a contradiction since we know that, by our assumptions, $\varphi(t,x)=Y^{t,x}_t$ a.e.. This conclude the proof.
\endproof

\begin{prop}\label{parsubopgensubvis}
The function $u(t,x)$ is a continuous viscosity subsolution to the \cref{parsubopgenprob}.
\end{prop}
\proof
We know from \cref{parsubopgenviscont} that $u$ is continuous, thus we just have to prove the subsolution property to conclude the proof.\\
Let $L$ be an element of $\mathcal L_F^t$ and $u_L$ as defined in \cref{parsemlinvis}, then if $\varphi$ is a supertangent to $u_L$ in $(t,x)$ we have that, by the definition of $\mathcal L_F^t$,
\[
F\left(t,x,\nabla_x\varphi,D^2_x\varphi\right)\ge L\left(t,x,\nabla_x\varphi,D^2_x\varphi\right)\ge-\partial_t\varphi(t,x)-f(t,x,u_L,\nabla_x\varphi),
\]
therefore $u_L$ is a viscosity subsolution to the \cref{parsubopgenprob} at $(t,x)$. Thanks to the arbitrariness of $t$, $L$ and $x$ we then have that $u_L$ is a viscosity subsolution in $(t,x)$ for any $L\in\mathcal L_F^t$, $x\in\mathds R^N$ and $t\in(0,T)$. From \cref{supeqrem} we have that
\[
\sup_{L\in\mathcal L_F^t}u_L(t,x)=\sup_{(b,\sigma)\in\mathcal A_F^t}\mathds E\left(Y^{t,x}_t\right)=\sup_{(b,\sigma)\in\mathcal A_F}\mathds E\left(Y^{t,x}_t\right)=u(t,x),
\]
therefore the family of functions $\{u_L\}$ is locally equibounded, thanks to \cref{parsubopgenviscont}. Well known properties of viscosity solutions hence yield that
\[
\sup_{L\in\mathcal L_F^t}u_L(t,x)=u(t,x)
\]
is a viscosity subsolution for any $(t,x)\in(0,T)\times\mathds R^N$.
\endproof

We conclude this \namecref{parprob} with our main statement.

\begin{theo}
The function $u(t,x):=\sup\limits_{(b,\sigma)\in\mathcal A_F}\mathds E\left(Y^{t,x}_t\right)$ is the only viscosity solution to the \cref{parsubopgenprob} satisfying polynomial growth condition such that $u(T,x)=g(x)$ for any $x$ in $\mathds R^N$.
\end{theo}
\proof
The uniqueness is a consequence of \cref{parsubopgencomp} and \eqref{eq:parvisbound}, hence we only have to show that $u$ is a viscosity solution.\\
From \cref{parsubopgensubvis} we know that $u$ is a continuous viscosity subsolution and it is easy to see that $u(T,x)=g(x)$ for any $x\in\mathds R^N$, so we only need to prove the supersolution property of $u$. Fixed $(t,x)$ in $(0,T)\times\mathds R^N$, let $\psi$ be a subtangent to $u$ in $(t,x)$ which we assume, without loss of generality, equal to $u$ at $(t,x)$ and $\delta$ a positive constant such that
\begin{equation}\label{eq:subtandis}
\psi(s,y)\le u(s,y)\qquad\text{for any }(s,y)\in[t,t+\delta]\times B_\delta(x).
\end{equation}
We know, thanks to \cref{charsublin}, that there exists a continuous and deterministic $L\in\mathcal L_F$ for which
\[
F\left(t,x,\nabla_x\psi,D^2_x\psi\right)=L\left(t,x,\nabla_x\psi,D^2_x\psi\right)
\]
and assume by contradiction
\begin{align*}
F\left(t,x,\nabla_x\psi,D^2_x\psi\right)=L\left(t,x,\nabla_x\psi,D^2_x\psi\right)>-\partial_t\psi(t,x)-f_\sigma(t,x,u,\nabla_x\psi\sigma).
\end{align*}
Then, by continuity,
\begin{equation}\label{eq:supsoldisass}
\partial_t\psi(s,y)+L\left(s,y,\nabla_x\psi,D^2_x\psi\right)>-f_\sigma(s,y,\psi,\nabla_x\psi\sigma)
\end{equation}
for any $(s,y)\in[t,t+\delta]\times B_\delta(x)$, possibly taking a smaller $\delta$.\\
We denote with $(b,\sigma)$ and $(X,Y,Z)$, respectively, the element of $\mathcal A$, which, to repeat, is continuous and deterministic, and the solution to the FBSDE \eqref{eq:fbsde} associated to $L$. We define the stopping time
\[
\tau:=(t+\delta)\wedge\inf\left\{s\in[t,\infty):\left|X_s^{t,x}-x\right|\ge\delta\right\},
\]
let $\left(\overline Y_s,\overline Z_s\right)$ be the solution to the BSDE
\[
\overline Y_s=u\left(\tau,X_\tau^{t,x}\right)+\int_{s\wedge\tau}^\tau f_\sigma\left(r,X^{t,x}_r,\overline Y_r,\overline Z_r\right)dr-\int_{s\wedge\tau}^\tau\overline Z_rdW_r,\quad s\in[t,T]
\]
and $\left(\hat Y_s,\hat Z_s\right):=\left(\psi\left(s,X^{t,x}_{s\wedge\tau}\right),(\nabla_x\psi\sigma)\left(s,X^{t,x}_{s\wedge\tau}\right)\right)$ which, by Itô's formula, is solution to
\[
\begin{aligned}
\hat Y_s=&\psi\left(\tau,X^{t,x}_\tau\right)-\int_{s\wedge\tau}^\tau\hat Z_rdW_r\\
&-\int_{s\wedge\tau}^\tau\left(\partial_t\psi\left(r,X^{t,x}_r\right)+L\left(r,X^{t,x}_r,\nabla_x\psi,D^2_x\psi\right)\right)dr,
\end{aligned}
\quad s\in[t,T].
\]
We know from the \cref{cauchydpp} that
\begin{equation}\label{eq:eqassum}
\sup_{(b,\sigma)\in\mathcal A_F}\mathds E\left(\overline Y_t\right)=u(t,x)=\psi(t,x),
\end{equation}
but by \eqref{eq:subtandis} we have $u\left(\tau,X^{t,x}_\tau\right)\ge\psi\left(\tau,X^{t,x}_\tau\right)$, which together with \eqref{eq:supsoldisass} imply, thanks to \cref{bsdeconf}, that $\overline Y_t>\psi(t,x)$ a.e., contradicting \eqref{eq:eqassum}.
\endproof

\section{Connection with 2BSDEs}\label{2bsdesec}

In this section we will briefly present a different approach to second order BSDEs, 2BSDEs for short, using our dynamic programming principle. This is intended as a short survey and not as a complete analysis of the subject.

We start giving the classical formulation of 2BSDE. Following \cite{sonertouzizhang11,possamai13} assume that $\Omega:=\left\{\omega\in C\left([0,T];\mathds R^N\right):\omega_0=0\right\}$ and let $\mathds P_0$ be the Wiener measure. Note that in this space the Brownian motion $W$ is a projection, i.e. $W_t(\omega)=\omega_t$. Denote with $[W]_t$ the quadratic variation of the projection and with
\[
\hat a_t:=\lim_{\varepsilon\downarrow0}\frac{[W]_t-[W]_{t-\varepsilon}}{\varepsilon}
\]
its variation. We will then denote with $\mathcal P_W$ the set of the probability measures $\mathds P$ such that $[W]$ is absolutely continuous in $t$ and $\hat a\in\mathds S^N_+$, $\mathds P$--a.e.. In particular $\mathds P_0\in\mathcal P_W$ because $[W]_t=tI_N$ and $\hat a_t=I_N$, $\mathds P_0$--a.e., where $I_N$ is the $N\times N$ identity matrix. Moreover, let $\mathcal P_S$ be the subset of $\mathcal P_W$ composed by the probability measures $\mathds P^\alpha:=\mathds P_0\circ(X^\alpha)^{-1}$, where
\[
X^\alpha_t:=\int_0^t\alpha^{1/2}_sdW_s,\qquad\mathds P_0\text{--a.e.},
\]
and $\alpha$ is a progressive process in $\mathds S^N_+$ such that, for two fixed $\underline a,\overline a$ in $\mathds S^N_+$, $\underline a\le\alpha\le\overline a$, $\mathds P_0$--a.e.. It is then apparent a link between $\mathcal P_S$ and the control set $\mathcal A_F$.\\
Now, given a function
\[
h:[0,T]\times\Omega\times\mathds R\times\mathds R^N\times D_h\to\mathds R,
\]
where $D_h$ is a subset of $\mathds R^{N\times N}$ containing 0, define for any $a\in\mathds S^N_+$
\[
f(t,\omega,y,z,a):=\sup_{\gamma\in D_h}\left(\frac12\langle\gamma,a\rangle-h(t,\omega,y,z,\gamma)\right).
\]
Furthermore let $\hat f(t,y,z):=f(t,y,z,\hat a_t)$,
\[
\mathcal P_h^2:=\left\{\mathds P\in\mathcal P_S:\mathds E^{\mathds P}\left(\int_0^T\left|\hat f(t,0,0)\right|^2dt\right)<\infty\right\}
\]
and $\mathcal P_h^2(t,\mathds P):=\left\{\mathds P'\in\mathcal P_h^2:\mathds P'=\mathds P\text{ on }\mathcal F_t\right\}$.

A pair of progressive processes $(Y,Z)$ is solution to the 2BSDE
\begin{equation}\label{eq:2bsde}
Y_t=\xi+\int_t^T\hat f(s,Y_s,Z_s)ds-\int^T_tZ_sdW_s+K_T-K_t
\end{equation}
if, for any $\mathds P\in\mathcal P_h^2$,
\begin{enumerate}[label=(\roman*)]
\item $Y_T=\xi$, $\mathds P$--a.e.;
\item the process $K^{\mathds P}$ defined below has non decreasing path $\mathds P$--a.e.,
\[
K_t^{\mathds P}=Y_0-Y_t-\int_0^t\hat f(s,Y_s,Z_s)ds+\int^t_0Z_sdW_s,\quad t\in[0,T],\,\mathds P\text{--a.e.};
\]
\item the family $\left\{K^{\mathds P},\mathds P\in\mathcal P^2_h\right\}$ satisfies the minimum condition
\begin{equation}\label{eq:kmincond}
K_t^{\mathds P}=\essi_{\mathds P'\in\mathcal P_h^2(t,\mathds P)}\mathds E^{\mathds P}\left(K_T^{\mathds P'}\middle|\mathcal F_t\right),\qquad\mathds P\text{--a.e.}.
\end{equation}
\end{enumerate}
Under suitable conditions the 2BSDE \eqref{eq:2bsde} admits a unique solution and, if we denote with $\left(Y^{\mathds P}(r,\xi),Z^{\mathds P}(r,\xi)\right)$ the solution to the BSDE
\[
Y_t=\xi+\int_t^r\hat f(s,Y_s,Z_s)ds-\int^r_tZ_sdW_s,\qquad t\in[0,r],\,\mathds P\text{--a.e.},
\]
it can be proved that $Y_t=\esss\limits_{\mathds P'\in\mathcal P_h^2(t,\mathds P)}Y_t^{\mathds P'}(r,Y_r)$ for any $r\in[t,T]$ and $\mathds P\in\mathcal P_h^2$. The last identity is a dynamic programming principle and can be seen as the connection between 2BSDE and our method.

Now we will show a different formulation of 2BSDEs, using controls instead of probability measures. Let $\mathcal A$ be our control set, made up by the progressive processes in $L^2([0,T]\times\Omega;B)$, where $B$ is a Banach space, and, for any $t\in[0,T]$,
\[
\mathcal A(t,\alpha):=\left\{\alpha'\in\mathcal A:\alpha'_s=\alpha_s\text{ for any }s\in[0,t]\right\}.
\]
Then define the function
\[
f:[0,T]\times\Omega\times\mathds R\times\mathds R^N\times B\to\mathds R
\]
and assume that there exists a $C>0$ such that
\[
|f(t,y,z,\alpha)-f(t,y,z,\alpha')|\le C(1+|y|+|z|)|\alpha-\alpha'|.
\]
We will also assume that, for each $\alpha\in\mathcal A$, $f_\alpha(t,y,z):=f(t,y,z,\alpha)$ satisfies \cref{bsdeassum} uniformly with respect to $\alpha$. We point out that these conditions are not intended to be minimal. We then have that the BSDEs
\[
Y^\alpha_t=\xi_\alpha+\int_t^Tf_\alpha(s,Y_s^\alpha,Z_s^\alpha)ds-\int_t^TZ^\alpha_sdW_s,\qquad t\in[0,T],
\]
admit a unique solution for any $\alpha\in\mathcal A$. If we moreover require that, for any $t\in[0,T]$, $\sup\limits_{\alpha\in\mathcal A}\mathds E\left(|Y^\alpha_t|^2\right)<\infty$ (which can be achieved if, for example, $B$ is compact), then by the \cref{cauchydpp},
\begin{equation}\label{eq:poordpp}
\overline Y_t^\alpha:=\esss_{\alpha'\in\mathcal A(t,\alpha)}Y^{\alpha'}_t=\esss_{\alpha'\in\mathcal A(t,\alpha)}Y^{\alpha'}_t\left(r,\overline Y_r^{\alpha'}\right),\quad\text{for any }0\le t\le r\le T,
\end{equation}
where $Y^{\alpha'}\left(r,\overline Y_r^{\alpha'}\right)$ is solution to the BSDE $\left(\overline Y^{\alpha'}_r,f_\alpha,r\right)$. It is easy to see that, for any $\alpha\in\mathcal A$, $\overline Y^\alpha$ is a continuous progressive process in $L^2$ and $\overline Y^\alpha_T=\xi_\alpha$ a.e..\\
Using the same arguments of \cite{possamai13} we that, for each $\alpha\in\mathcal A$, there exist two progressive processes in $L^2$, $Z^\alpha$ and $K^\alpha$, such that $K^\alpha$ is a continuous and increasing process in $t$ with $K_0^\alpha=0$ and
\[
\overline Y_t^\alpha=\xi_\alpha+\int_t^Tf_\alpha\left(s,\overline Y^\alpha_s,\overline Z_s^\alpha\right)ds-\int^T_t\overline Z_s^\alpha dW_s+K^\alpha_T-K^\alpha_t,\qquad t\in[0,T].
\]
We also have that, as in \eqref{eq:kmincond},
\[
K_t^\alpha=\essi_{\alpha'\in\mathcal A(t,\alpha)}\mathds E\left(K_T^{\alpha'}\middle|\mathcal F_t\right),\qquad\text{a.e. for any }t\in[0,T],\,\alpha\in\mathcal A.
\]

If we let $X$ be as in \eqref{eq:fbsde} and $\left(\overline Y,\overline Z,K\right)$ be the solution to the 2BSDE (we omit the dependence on the control set $\mathcal A_F$)
\[
\overline Y_s^{t,x}=g\left(X^{t,x}_T\right)+\int_s^Tf_\sigma\left(r,X^{t,x}_T,\overline Y_r^{t,x},\overline Z_r^{t,x}\right)ds-\int^T_s\overline Z_r^{t,x}dW_r+K^{t,x}_T-K^{t,x}_s
\]
for any $x\in\mathds R^N$ and $0\le t\le s\le T$, we then have that the function $u(t,x):=\overline Y^{t,x}_t$ is the viscosity solution to \cref{parsubopgenprob}.

\begin{appendices}

\section{Some Probability Results}

Here we give some probability results we use in this paper.

\begin{lem}\label{cadlagcont}
Let $\{U_t\}_{t\in[0,\infty)}$ be a cadlag process, then for any $\varepsilon>0$ there exists a $\delta>0$ such that
\[
\mathds P(\{|U_t-U_s|<\varepsilon,\text{ for any }s\in[t,t+\delta)\})>0.
\]
\end{lem}
\proof
Our argument is by contradiction. Assume that there exists an $\varepsilon>0$ such that for any $\delta>0$
\[
\mathds P(\{|U_t-U_s|<\varepsilon,\text{ for any }s\in[t,t+\delta)\})=0,
\]
which is equivalent to
\[
\mathds P(\{|U_t-U_s|\ge\varepsilon,\text{ for any }s\in[t,t+\delta)\})=1.
\]
Let, for any positive integer $n$,
\[
A_n:=\left\{|U_t-U_s|\ge\varepsilon,\text{ for any }s\in\left[t,t+\frac1n\right)\right\},
\]
then $A_{n}\subseteq A_{k}$ if $k\le n$ and
\[
A:=\bigcap_{n=1}^\infty A_n=\left\{\lim_{s\downarrow t}|U_t-U_s|\ge\varepsilon\right\}.
\]
Since $U$ is right continuous we know that $\mathds P(A)=0$ which contradicts our assumption, since $\mathds P(A)=\lim\limits_{n\to\infty}\mathds P(A_n)=1$.
\endproof

Consider the SDE
\begin{equation}\label{eq:sde}
X^{t,\zeta}_s=\zeta+\int_t^s\sigma\left(r,X^{t,\zeta}_r\right)dW_r+\int_t^sb\left(r,X^{t,\zeta}_r\right)dr,\qquad s\in[t,\infty),
\end{equation}
under the following assumptions:
\begin{assum}\label{sdeassum}
$t\in[0,\infty)$, $\zeta\in L^2\left(\Omega,\mathcal F_t;\mathds R^N\right)$ and for the functions
\[
b:[0,\infty)\times\Omega\times\mathds R^N\to\mathds R^N\text{ and }\sigma:[0,\infty)\times \Omega\times\mathds R^N\to\mathds R^{N\times M}
\]
there exists a positive constant $\ell$ such that a.e., for any $r\in[0,\infty)$, $x,y\in\mathds R^N$,
\begin{enumerate}[label=(\roman*)]
\item $|b(r,x)-b(r,y)|+|\sigma(r,x)-\sigma(r,y)|\le\ell|x-y|$;
\item $\{(b,\sigma)(t,0)\}_{t\in[0,\infty)}$ is a progressive process belonging to $L^2([0,T]\times\Omega)$ for any $T\in[0,\infty)$.
\end{enumerate}
\end{assum}

A solution to this SDE is a continuous progressive process $X$ as in \eqref{eq:sde} such that $X\in L^2([0,T]\times\Omega)$ for any $T\in[0,\infty)$. The next theorem summarizes some SDE result given in \cite{krylovbook}.

\begin{theo}\label{sdeex}
Under \cref{sdeassum} there exists a unique solution to the SDE \eqref{eq:sde}. Moreover, for any $T\in[t,\infty)$, there exists a constant $c$, depending only on $\ell$ and $T$ such that
\begin{gather*}
\mathds E\left(\sup_{s\in[t,T]}\left|X^{t,\zeta}_s-\zeta\right|^2\right)\le c\mathds E\left(\int_t^T\left(|b(s,0)|^2+|\sigma(s,0)|^2\right)ds\right),\\
\mathds E\left(\sup_{s\in[t,T]}\left|X^{t,\zeta}_s\right|^2\right)\le c\mathds E\left(|\zeta|^2+\int_t^T\left(|b(s,0)|^2+|\sigma(s,0)|^2\right)ds\right),\\
\mathds E\left(\sup_{s\in[t,T]}\left|X^{t,\zeta}_s-X^{t,\zeta'}_s\right|^2\right)\le c\mathds E\left(|\zeta-\zeta'|^2\right).
\end{gather*}
\end{theo}

\begin{rem}\label{sdestopstart}
The results obtained in this section hold even for SDEs with an a.e. finite stopping time $\tau$ as starting time. In fact if for any $\zeta$ in $L^2\left(\Omega,\mathcal F_\tau;\mathds R^N\right)$ we define
\[
\overline b(t,x):=b(t,x+\zeta)\chi_{\{\tau\le t\}}\text{ and }\,\overline\sigma(t,x):=\sigma(t,x+\zeta)\chi_{\{\tau\le t\}},
\]
then $X^{\tau,\zeta}$ is solution of the SDE $(b,\sigma)$ if and only if $\overline X^{0,0}:=X^{\tau,\zeta}-\zeta$ is solution of the SDE $\left(\overline b,\overline\sigma\right)$. The claim can be easily obtained from this.
\end{rem}

\begin{rem}
By the strong Markov property, for any a.e. finite stopping time $\tau$, the process $\{W_t^\tau\}_{t\in[0,\infty)}:=\{W_{\tau+t}-W_\tau\}_{t\in[0,\infty)}$ is a Brownian motion. Thus if $b$ and $\sigma$ are are progressive with respect to the filtration $\{\mathcal F^\tau_t\}_{t\in[0,\infty)}$ then any solution to the SDE $(b,\sigma)$ with initial data $\tau+t$ and $\zeta\in L^2\left(\Omega,\mathcal F_t^\tau;\mathds R^N\right)$ is also progressive with respect to that filtration. In fact, in this case, the stochastic integral with respect to $W_t^\tau$ is the same as the one with respect to $W_{\tau+t}$.
\end{rem}

\subsection{Backward Stochastic Differential Equations}

In this subsection we give some results on BSDEs used in our investigation. Most of them are well known and actually hold under more general assumptions. We refer to \cite{syspardoux,phamsdebook,pardouxbook,brianddelyonhupardouxstoica03} for their proofs.

We will work under the followings assumptions:

\begin{assum}\label{bsdeassum}
Let $T\in[0,\infty)$, $\xi\in L^2\left(\Omega,\mathcal F_T,\mathds P;\mathds R^M\right)$ and
\[
f:[0,T]\times \Omega\times\mathds R^M\times\mathds R^{M\times N}\to\mathds R^M
\]
a function which admits a positive constant $\ell$ and a real number $\mu$ such that a.e., for any $t\in[0,T]$, $y,y'\in\mathds R^M$ and $z,z'\in\mathds R^{M\times N}$,
\begin{enumerate}[label=(\roman*)]
\item $\{f(s,0,0)\}_{s\in[0,T]}$ is a progressive process belonging to $L^2([0,T]\times\Omega)$;
\item $|f(t,y,z)|\le|f(t,0,z)|+\ell(1+|y|)$;
\item $|f(t,y,z)-f(t,y,z')|\le\ell|z-z'|$;
\item $(y-y')(f(t,y,z)-f(t,y',z))\le\mu|y-y'|^2$;
\item $v\mapsto f(t,v,z)$ is continuous.
\end{enumerate}
\end{assum}

A solution to the BSDE $(\xi,f,T)$, where $\xi$ and $T$ have respectively the role of a final condition and a terminal time, is a pair $(Y,Z)$ of progressive processes belonging to $L^2([0,T]\times\Omega)$ such that
\begin{equation}\label{eq:bsde}
Y_t=\xi+\int_t^Tf(s,Y_s,Z_s)ds-\int_t^TZ_sdW_s,\quad\text{for any }t\in[0,T].
\end{equation}

The followings are classical results of BSDE theory.

\begin{theo}\label{bsdeex}
Under the \cref{bsdeassum} the BSDE \eqref{eq:bsde} has a unique solution $(Y,Z)$. Furthermore there exists a constant $c$, which depends on $T$, $\mu$ and $\ell$, such that
\[
\mathds E\left(\sup_{t\in[0,T]}|Y_t|^2+\int_0^T|Z_t|^2dt\right)\le c\mathds E\left(|\xi|^2+\int_0^T|f(t,0,0)|^2dt\right)
\]
and, if $(Y',Z')$ is the solution to the BSDE $(\xi',f',T)$,
\begin{multline*}
\mathds E\left(\sup_{t\in[0,T]}|Y_t-Y'_t|^2+\int_0^T|Z_t-Z'_t|^2dt\right)\\
\le c\mathds E\left(|\xi-\xi'|^2+\int_0^T|f(t,Y'_t,Z'_t)-f'(t,Y'_t,Z'_t)|^2dt\right).
\end{multline*}
\end{theo}
\begin{theo}\label{bsdeconf}
Assuming $M=1$, let $(Y,Z)$ be the solution to the BSDE $(\xi,f,T)$ under the \cref{bsdeassum} and
\[
Y'_t=\xi'+\int_t^TV_sds-\int_t^TZ'_sdW_s,\qquad t\in[0,T],
\]
where $\xi'\in L^2\left(\Omega,\mathcal F_T,\mathds P;\mathds R\right)$, $Y',V\in L^2([0,T]\times\Omega)$ and $Z'\in L^2([0,T]\times\Omega)$. Suppose that $\xi\le\xi'$ a.e. and $f(t,Y'_t,Z'_t)\le V_t$ a.e. for the $dt\times d\mathds P$ measure. Then, for any $t\in[0,T]$, $Y_t\le Y'_t$ a.e..\\
If moreover $Y_0=Y'_0$ a.e., then $Y_t=Y'_t$ a.e. for any $t\in[0,T]$, or in other words, whenever either $\mathds P(\{\xi<\xi'\})>0$ or $f(s,Y'_s,Z'_s)<V_s$, for any $(y,z)$ in $\mathds R\times\mathds R^N$ on a set of positive $dt\times d\mathds P$ measure, then $Y_0<Y'_0$.
\end{theo}
\begin{prop}\label{stopbsde}
Let $(Y,Z)$ be the solution to the BSDE \eqref{eq:bsde} and assume that there exists a stopping time $\tau$ such that $\tau\le T$, $\xi$ is $\mathcal F_\tau$--measurable and $f(t,y,z)=0$ on the set $\{\tau\le t\}$. Then $Y_t=Y_{\tau\wedge t}$ a.e. and $Z_t=0$ a.e. on the set $\{\tau\le t\}$.
\end{prop}

\section{Comparison Theorem}

Consider the parabolic problem
\begin{equation}\label{eq:defparpde}
\partial_tu(t,x)+F\left(t,x,u,\nabla u,D^2u\right)=0,\quad t\in(0,T),x\in\mathds R^N,
\end{equation}
where $F$ is a continuous elliptic operator which admits, for any $t\in[0,T]$, $(x,r,p,S)$ and $(y,r',p',S')$ in $\mathds R^N\times\mathds R\times\mathds R^N\times\mathds S^N$, a $\mu\in\mathds R$ and a positive constant $\ell$ such that
\begin{enumerate}[label=(\roman*)]
\item $|F(t,x,r,p,S)-F(t,x,r,p,S')|\le\ell\left(1+|x|^2\right)|S-S'|$;
\item $|F(t,x,r,p,S)-F(t,x,r,p',S)|\le\ell(1+|x|)|p-p'|$;
\item $|F(t,x,r,p,S)-F(t,y,r,p,S)|\le\ell(1+|x|+|y|)|x-y|\|(p,S)\|$;
\item\emph{(Monotonicity)} $(F(t,x,r,p,S)-F(t,x,r',p,S))(r-r')\le\mu|r-r'|^2$;
\item\label{domteocond} the continuity of the function $r\mapsto F(t,x,r,p,S)$ is independent from the fourth variable.
\end{enumerate}
Notice that, given a compact set $K\subset\mathds R^N\times\mathds R\times\mathds S^N$, \cref{domteocond} and the Heine--Cantor theorem yield the existence of a modulus of continuity $\omega_K$ such that, if $(x,r,S),(x,r',S)\in K$,
\begin{equation}\label{eq:domteocond}
|F(t,x,r,p,S)-F(t,x,r',p,S)|\le\omega_K(|r-r'|)
\end{equation}
for any $t\in[0,T]$ and $p\in\mathds R^N$.

Here we give a comparison result, which is an adaptation of \cite[Theorem C.2.3]{book:peng} for problem \eqref{eq:defparpde}. To prove it we will need the following \namecref{domteo}, adaptation of \cite[Theorem C.2.2]{book:peng} which can be proved similarly. Notice that condition \textbf{(G)} in \cite{book:peng} is replaced by \eqref{eq:domteo.1}.

\begin{theo}\label{domteo}
Let $\{F_i\}_{i=1}^k$ be a collection of continuous functions from $[0,T]\times\mathds R^N\times\mathds R\times\mathds R^N\times\mathds S^N$ to $\mathds R^N$ and assume that, if $(x,r,S),(y,r,S)$ belong to a compact set $K\subset\mathds R^N\times\mathds R\times\mathds S^N$, there exist a constant $C_K$ and a modulus of continuity $\omega_K$ such that, for any $i\in\{1,\cdots,k\}$,
\begin{equation}\label{eq:domteo.1}
|F_i(t,x,r,p,S)-F_i(t,y,r,p,S)|\le C_K(1+|p|)|x-y|+\omega_K(|x-y|).
\end{equation}
Furthermore assume the following domination condition: there exists a collection of positive constants $\{\beta_i\}_{i=1}^k$ satisfying
\[
\sum_{i=1}^k\beta_iF_i(t,x,r_i,p_i,S_i)\le0
\]
for each $(t,x)\in[0,T]\times\mathds R^N$ and $(r_i,p_i,S_i)$ such that $\sum\limits_{i=1}^k\beta_ir_i\ge0$, $\sum\limits_{i=1}^k\beta_ip_i=0$ and $\sum\limits_{i=1}^k\beta_iX_i\le0$.\\
For any $i\in\{1,\cdots,k\}$, let $u_i$ be a viscosity subsolution of
\[
\partial_tu(t,x)+F_i\left(t,x,u,\nabla u,D^2u\right)=0,\quad t\in(0,T),x\in\mathds R^N,
\]
and assume that $\sum\limits_{i=1}^k\beta_iu_i(T,\cdot)\le0$ and $\left(\sum\limits_{i=1}^k\beta_iu_i(\cdot,x)\right)^+\to0$ uniformly as $|x|\to\infty$. Then $\sum\limits_{i=1}^k\beta_iu_i(t,\cdot)\le0$ for any $t\in(0,T)$.
\end{theo}

\begin{theo}\label{defviscomp}
Let $u$ and $v$ be, respectively, a viscosity subsolution and a viscosity supersolution to \eqref{eq:defparpde} satisfying polynomial growth condition. Then, if $u|_{t=T}\le v|_{t=T}$, $u\le v$ on $(0,T]\times\mathds R^N$.
\end{theo}
\proof
We set $\phi(x):=(1+|x|^2)^\frac c2$,
\[
\gamma>\mu+\ell\sup_{x\in\mathds R^N}\left((1+|x|)\frac{|\nabla\phi(x)|}{\phi(x)}+\left(1+|x|^2\right)\frac{|D^2\phi(x)|}{\phi(x)}\right),
\]
$u_1(t,x):=\dfrac{e^{-\gamma t}u(t,x)}{\phi(x)}$ and $u_2(t,x):=-\dfrac{e^{-\gamma t}v(t,x)}{\phi(x)}$, where $c$ is such that both $|u_1|$ and $|u_2|$ converge uniformly to 0 as $|x|\to\infty$. Notice that
\[
\nabla\phi(x)=c\frac{\phi(x)x^\dag}{1+|x|^2}\quad\text{and}\quad D^2\phi(x)=\phi(x)\left(\frac c{1+|x|^2}I+\frac{c(c-2)}{(1+|x|^2)^2}x\otimes x\right),
\]
therefore $\gamma$ is well defined. We also set the operators $F_1(t,x,r,p,S)$, given by
\[
\frac{e^{-\gamma t}}{\phi(x)}F\!\left(t,x,e^{\gamma t}r\phi,e^{\gamma t}\left(r\nabla\phi+\phi p^\dag\right),e^{\gamma t}\left(rD^2\phi+\nabla\phi\otimes p+p\otimes\nabla\phi+\phi S\right)\right)\!,
\]
and $F_2(t,x,r,p,S)$, given by
\[
-\frac{e^{-\gamma t}\!}{\phi(x)}F\!\left(\!t,x,\!-e^{\gamma t}r\phi,\!-e^{\gamma t}\!\left(r\nabla\phi\!+\!\phi p^\dag\right)\!,\!-e^{\gamma t}\!\left(rD^2\phi\!+\!\!\nabla\phi\!\otimes\!p+p\!\otimes\!\nabla\phi\!+\!\phi S\right)\right)\!.
\]
It is easy to check that, for $i\in\{1,2\}$, $F_i$ is still continuous, elliptic, Lipschitz continuous in $p$ and $S$, that its monotonicity constant is
\[
\mu+\ell\sup_{x\in\mathds R^N}\left((1+|x|)\frac{|\nabla\phi(x)|}{\phi(x)}+\left(1+|x|^2\right)\frac{|D^2\phi(x)|}{\phi(x)}\right),
\]
i.e. is lower than $\gamma$, and $u_i$ is a viscosity subsolution to
\[
\partial_tu(t,x)-\gamma u(t,x)+F_i\left(t,x,u,\nabla u,D^2u\right)=0,\quad t\in(0,T),x\in\mathds R^N.
\]
It can also be checked that, if $(x,r,S),(y,r,S)$ belong to a compact set $K\subset\mathds R^N\times\mathds R\times\mathds S^N$, there exist a constant $C_K$ and a modulus of continuity $\widetilde\omega_K$ bigger than $\omega_K$ in \eqref{eq:domteocond} such that
\[
|F_i(t,x,r,p,S)-F_i(t,y,r,p,S)|\le C_K(1+|p|)|x-y|+\widetilde\omega_K(|x-y|)
\]
for any $i\in{1,2}$, $t\in[0,T]$ and $p\in\mathds R^N$.\\
Furthermore $(u_1+u_2)|_{t=T}\le0$ and $F_1(t,x,r,p,S)+F_2(t,x,-r,-p,-S)=0$. From these properties we have, for any $r_i\in\mathds R$, $p_i\in\mathds R^N$ and $S_i\in\mathds S^N$ such that $r=r_1+r_2\ge0$, $p_1=-p_2$ and $S=S_1+S_2\le0$,
\begin{multline*}
-\gamma r+F_1(t,x,r_1,p_1,S_1)+F_2(t,x,r_2,p_2,S_2)\\
\begin{aligned}
=&-\gamma r+F_1(t,x,r_1,p_1,S_1)+F_2(t,x,-r_1,-p_1,-S_1)\\
&+(F_2(t,x,r_2,p_2,S_2)-F_2(t,x,r_2-r,p_2,S_2-S))\frac rr\\
\le&-\gamma r+(F_2(t,x,r_2,p_2,S_2)-F_2(t,x,r_2-r,p_2,S_2))\frac rr\\
\le&-\gamma r+\gamma r=0.
\end{aligned}
\end{multline*}
As a consequence we have that all the conditions of \cref{domteo} are satisfied, thus $u_1+u_2\le0$, or equivalently, $u\le v$ in $(0,T]\times\mathds R^N$.
\endproof

\end{appendices}

\phantomsection

\pdfbookmark[1]{References}{References}
\emergencystretch=2em
\bibliography{biblio}

\begin{thebibliography}{10}

\bibitem{bellmanmatrixbook}
R.~Bellman.
\newblock {\em Introduction to matrix analysis}.
\newblock Society for Industrial and Applied Mathematics, Philadelphia, 1997.

\bibitem{bensoussa82}
A.~Bensoussan.
\newblock {\em Lectures on stochastic control}, pages 1--62.
\newblock Springer Berlin Heidelberg, 1982.

\bibitem{brianddelyonhupardouxstoica03}
P.~Briand, B.~Delyon, Y.~Hu, E.~Pardoux, and L.~Stoica.
\newblock Lp solutions of backward stochastic differential equations.
\newblock {\em Stochastic Processes and their Applications}, 108(1):109--129,
  11 2003.

\bibitem{int2bsde}
P.~Cheridito, H.~M. Soner, N.~Touzi, and N.~Victoir.
\newblock Second order backward stochastic differential equations and fully
  nonlinear parabolic pdes.
\newblock {\em Communications on Pure and Applied Mathematics},
  60(7):1081--1110, 2007.

\bibitem{userguide}
M.~G. Crandall, H.~Ishii, and P.-L. Lions.
\newblock User{\rq}s guide to viscosity solutions of second order partial
  differential equations.
\newblock {\em Bulletin of the American Mathematical Society}, 27(1):1--68, 11
  1992.

\bibitem{Delarue_2002}
F.~Delarue.
\newblock On the existence and uniqueness of solutions to fbsdes in a
  non--degenerate case.
\newblock {\em Stochastic Processes and their Applications}, 99(2):209--286, 6
  2002.

\bibitem{Delarue_Guatteri_2006}
F.~Delarue and G.~Guatteri.
\newblock Weak existence and uniqueness for forward--backward sdes.
\newblock {\em Stochastic Processes and their Applications},
  116(12):1712--1742, 12 2006.

\bibitem{art:pathpeng}
L.~Denis, M.~Hu, and S.~Peng.
\newblock Function spaces and capacity related to a sublinear expectation:
  Application to g-brownian motion paths.
\newblock {\em Potential Analysis}, 34(2):139--161, 5 2010.

\bibitem{ekrentouzizhang16}
I.~Ekren, N.~Touzi, and J.~Zhang.
\newblock Viscosity solutions of fully nonlinear parabolic path dependent pdes:
  Part i.
\newblock {\em The Annals of Probability}, 44(2):1212--1253, 3 2016.

\bibitem{ekrentouzizhang16_2}
I.~Ekren, N.~Touzi, and J.~Zhang.
\newblock Viscosity solutions of fully nonlinear parabolic path dependent pdes:
  Part ii.
\newblock {\em The Annals of Probability}, 44(4):2507--2553, 7 2016.

\bibitem{flemingrishel75}
W.~Fleming and R.~Rishel.
\newblock {\em Deterministic and Stochastic Optimal Control}.
\newblock Springer New York, New York, NY, 1975.

\bibitem{kazi-tanipossamaizhou15}
N.~Kazi-Tani, D.~Possamaï, and C.~Zhou.
\newblock Second-order bsdes with jumps: Formulation and uniqueness.
\newblock {\em The Annals of Applied Probability}, 25(5):2867--2908, 10 2015.

\bibitem{krylovbook}
N.~V. Krylov.
\newblock {\em Controlled diffusion processes}.
\newblock Springer, New York, 2009.

\bibitem{hjlions2}
P.~L. Lions.
\newblock Optimal control of diffusion processes and hamilton--jacobi--bellman
  equations part 2: viscosity solutions and uniqueness.
\newblock {\em Communications in Partial Differential Equations},
  8(11):1229--1276, 1 1983.

\bibitem{hjlions1}
P.~L. Lions.
\newblock Optimal control of diffustion processes and hamilton-jacobi-bellman
  equations part i: the dynamic programming principle and application.
\newblock {\em Communications in Partial Differential Equations},
  8(10):1101--1174, 1 1983.

\bibitem{Ma_Yong_2007}
J.~Ma and J.~Yong.
\newblock {\em Forward--Backward Stochastic Differential Equations and their
  Applications}.
\newblock Springer Berlin Heidelberg, 2007.

\bibitem{matoussipossamaizhou13}
A.~Matoussi, D.~Possamaï, and C.~Zhou.
\newblock Robust utility maximization in nondominated models with 2bsde: the
  uncertain volatility model.
\newblock {\em Mathematical Finance}, 25(2):258--287, 6 2013.

\bibitem{nisio15}
M.~Nisio.
\newblock {\em Stochastic Control Theory}.
\newblock Springer Japan, 2015.

\bibitem{syspardoux}
{\'E}.~Pardoux.
\newblock {\em Backward Stochastic Differential Equations and Viscosity
  Solutions of Systems of Semilinear Parabolic and Elliptic PDEs of Second
  Order}, pages 79--127.
\newblock Birkh{\"a}user Boston, 1998.

\bibitem{Pardoux_Peng_1990}
{\'E}.~Pardoux and S.~Peng.
\newblock Adapted solution of a backward stochastic differential equation.
\newblock {\em Systems \& Control Letters}, 14(1):55--61, 1 1990.

\bibitem{pardouxbook}
{\'E}.~Pardoux and A.~R{\u a}{\c s}canu.
\newblock {\em Stochastic Differential Equations, Backward SDEs, Partial
  Differential Equations}.
\newblock Springer International Publishing, 2014.

\bibitem{math/0601035}
S.~Peng.
\newblock {\em G-Expectation, G-Brownian Motion and Related Stochastic Calculus
  of It{\^o} Type}, pages 541--567.
\newblock Springer Berlin Heidelberg, 2007.

\bibitem{book:peng}
S.~Peng.
\newblock {\em Nonlinear Expectations and Stochastic Calculus under
  Uncertainty}.
\newblock Springer Berlin Heidelberg, 2019.

\bibitem{phamsdebook}
H.~Pham.
\newblock {\em Continuous-time Stochastic Control and Optimization with
  Financial Applications}.
\newblock Springer, Berlin Heidelberg, 2009.

\bibitem{possamai13}
D.~Possamaï.
\newblock Second order backward stochastic differential equations under a
  monotonicity condition.
\newblock {\em Stochastic Processes and their Applications}, 123(5):1521--1545,
  5 2013.

\bibitem{possamaitanzhou18}
D.~Possamaï, X.~Tan, and C.~Zhou.
\newblock Stochastic control for a class of nonlinear kernels and applications.
\newblock {\em The Annals of Probability}, 46(1):551--603, 1 2018.

\bibitem{sqmatrix}
B.~A. Schmitt.
\newblock Perturbation bounds for matrix square roots and pythagorean sums.
\newblock {\em Linear Algebra and its Applications}, 174:215--227, 9 1992.

\bibitem{schneiderconvex}
R.~Schneider.
\newblock {\em Convex bodies: the Brunn-Minkowski theory}.
\newblock Cambridge University Press, 2014.

\bibitem{sonertouzizhang11}
H.~M. Soner, N.~Touzi, and J.~Zhang.
\newblock Wellposedness of second order backward sdes.
\newblock {\em Probability Theory and Related Fields}, 153(1--2):149--190, 2
  2011.

\bibitem{williams08}
D.~Williams.
\newblock {\em Probability with martingales}.
\newblock Cambridge University Press, Cambridge New York, 2008.

\end{thebibliography}
\end{document}